\RequirePackage{fix-cm}
\documentclass[smallextended]{svjour3}       
\smartqed  
\usepackage{graphicx}
\usepackage{amssymb}
\usepackage{amsmath}
\begin{document}

\title{The Subdifferential Descent Method \\in a Nonsmooth Variational Problem \thanks{}
}

\titlerunning{The Subdifferential Descent Method in a Nonsmooth Variational Problem}        

\author{A. V. Fominyh }


\institute{A. V. Fominyh \at
              St. Petersburg State University, 7/9 Universitetskaya nab., St. Petersburg, 199034, Russia \\
              Tel.: +7 905 212 60 65\\
              \email{alexfomster@mail.ru} }          

\date{Received: date / Accepted: date}

\maketitle

\begin{abstract}
The paper is devoted to the classical variational problem with a nonsmooth integrand of the functional to be minimized. The integrand is supposed to be subdifferentiable. Under some natural conditions the subdifferentiability of the functional considered is proved. The problem of finding the subdifferential descent is being solved and the subdifferential descent method is applied to solve the original problem. The algorithm developed is demonstrated by examples. 

\end{abstract}

\section{Introduction}

Most of existing numerical methods for solving problems of the calculus of variations are developed for the case when the integrand is continuously differentiable with respect to the variables sought. This paper is aimed at solving the simplest variational problem under the assumption that the integrand of the minimized functional is nonsmooth and only subdifferentiable with respect to the unknown function and to its derivative. 

To study the existence of generalized Bolza problem solution (whose particular case is the problem considered in this paper), a deep theory is constructed in papers \cite{Rock1}, \cite{Rock2}, \cite{Rock3}, which uses Fenchel-Moreau duality and other interesting facts from convex analysis, as well as some nontrivial results of functional analysis. The proof for generalized Bolza problem in the case of delay can be found in \cite{OrtizWolenski}. Some general results on solution existence are also contained in paper \cite{Ioffe1}. 

Most of the works which consider nonsmooth variational problems are of a theoretical nature and investigate necessary and sufficient conditions for a minimum. For example, in paper \cite{IoffeRockafellar} the necessary conditions are formulated in terms of subgradients and generalize the well-known Euler and Weirstrass conditions of the classical theory of the variational calculus. In paper \cite{LoewenRockafellar} the subdifferential of the Hamiltonian is used to formulate the necessary minimum conditions, and the Hamiltonian is a function conjugate to the Lagrangian, that is the duality theory is also used here. The necessary conditions (obtained in a similar form, as well as in Erdman's form) are investigated by qualitatively different methods in work \cite{Clarke}. The general necessary optimality conditions for the generalized Bolza problem were obtained in terms of the special differential constructions in book \cite{Morduhovich1} and for nonconvex differential inclusions --- in paper \cite{Morduhovich2}. In paper \cite{Zeidan} some sufficient conditions for a minimum are constructed based on the strengthened Weierstrass conditions. In work \cite{Dolgopolik1} the minimum conditions for various nonsmooth variational problems are obtained in terms of codifferentials. The results on the necessary minimum conditions for both generalized Bolza problem and the control problem of a differential inclusion are contained in recent work \cite{Ioffe2} in a fairly general and complete form; and it is interesting how the proof for the second problem is carried out by reducing it to the first one by adding an integral term of a special structure. In recent paper \cite{DolgopolikNew} minimum conditions have been obtained for variational problems with isoperimetric constraints; and these conditions are (in some cases) stronger than those known earlier. 

Paper \cite{Teo1} considered a special case of nonsmooth variational problems with equality and inequality constraints on integrand variables. Paper \cite{Teo2} also studied optimal control problems subject to nonsmooth functional constraints. In both of these papers some kind of smoothing technique was used in order to construct a numerical method for solving these problems. In works \cite{Tamasyan1}, \cite{TamasyanDemyanov} the methods of the subdifferential and the hypodifferential descents were applied to some classes of smooth variational problems with nonsmooth penalty summands which take into account the restriction on the right endpoint. These methods were also applied to constructing optimal control in problems with the subdifferentiable quality functional in paper \cite{Fominyh1} and also to the problem of transferring a system of differential equations from one point to another in works \cite{Fominyh3}, \cite{Fominyh4}. The finite-dimensional quasidifferential descent method was applied to optimization of a control system with a nonsmooth objective functional in Mayer form in paper \cite{Fominyh5}.  Despite the fact that in the last works listed the quality functional is subdifferentiable, it has a special structure (for example, being the maximum of Gateaux differentiable functionals); therefore, the calculation of its subdifferential is quite simple. In this paper the integrand of the functional to be minimized is nondifferentiable; therefore, the technique of the described papers is not applicable in this case. The key idea in overcoming this difficulty and obtaining a subdifferential in a constructive form (`` constructiveness '' here means the possibility of constructing an algorithm for solving a given problem) is to consider the trajectory and its derivative as independent variables (since, in fact, these variables are, of course, related to each other, we construct a penalty function of a special kind in order to take this relation into account (see the section Reduction to an Unconstrained Minimization Problem)).

\section{Statement of the Problem}  

In the paper we will use the following notations. $C_{n} [0, T]$ is a space of $n$-dimensional continuous on $[0, T]$ vector-functions, which are piecewise continuously differentiable with bounded on its domain derivative; $P_{n} [0, T]$ is a space of piecewise continuous and bounded on $[0, T]$ $n$-dimensional vector-functions. Denote $L_p^n [0, T]$, $1 \leqslant p < \infty$, the space of measurable on $[0, T]$ $n$-dimensional vector-functions which are $p$-summable and $L_\infty^n [0, T]$ --- the space of measurable on $[0, T]$ and almost everywhere bounded $n$-dimensional vector-functions. Denote $ \mathrm{co} P $ the convex hull of the set $P$. Let $B_r(c)$ ($D_r(c)$) denote a closed (open) ball in corresponding space with the radius $r$ and the center $c$; for some set $C$ in this space $B_r(C)$ ($D_r(C)$) denotes the union of all closed (open) balls with the radius $r$ and the centers from the set $C$. Denote $\langle a, b \rangle$ the scalar product of the vectors $a$, $b$ $\in R^d$. Let $X$ be a normed space, then $||\cdot||_X$ denotes the norm in this space and $X^*$ denotes the space conjugate to the space $X$. Finally, for some number $\alpha \in R$ let $o(\alpha)$ denote such a value that $o (\alpha) / \alpha \rightarrow 0 \ \text{if} \ \alpha \rightarrow 0$.

Let $x(t)$ be a piecewise continuously differentiable vector-function. Let \linebreak $t_0 \in [0, T)$ be a point of nondifferentiability of the vector-function $x(t)$, then for definiteness we assume that $\dot x(t_0)$ is a right-hand derivative of the vector-function $x(t)$ at the point~$t_0$. Similarly, we assume that $\dot x(T)$ is a left-hand derivative of the vector-function $x(t)$ at the point~$T$. As the derivative $\dot x(t)$ is supposed to be bounded on its domain, by previous paragraph notation we can assume that the vector-function $\dot x(t)$ belongs to the space $P_{n} [0, T]$.

Consider the following variational problem: it is required to minimize the functional
\begin{equation}
\label{1}
\overline J(x) = \int_0^T f(x(t), \dot x(t), t) dt
\end{equation}
with the boundary constraints
\begin{equation}
\label{2}
x(0) = x_{0}, \quad x(T) = x_T.
\end{equation}
In formula (\ref{1}) $f(x,\dot x,t)$, $t \in [0, T]$, is a given function, $T > 0$ is a given finite moment of time, $x(t)$ is an $n$-dimensional continuous vector-function, which is continuously differentiable at each $t \in [0, T]$ with the exception, possibly, of the finite number of points, and we suppose that its derivative is bounded on its domain. The function $f(x,\dot x,t)$ is continuous in $(x,\dot x,t)$ and locally Lipschitz continuous in $(x, \dot x)$ at each fixed point $t \in [0, T]$. In formula (\ref{2}) $x_0, x_T \in R^n$ are given vectors. 

In this paper we use both subdifferentials of functions in a finite-dimensional space and subdifferentials of functionals in a functional space. Despite the fact that the second concept generalizes the first one, for convenience we separately introduce definitions for both of these cases and for those specific functions (functionals) and their variables and spaces which are considered in the paper. 

Consider the space $R^n \times R^n$ with the standard norm. Let \linebreak $g = [g_1, g_2] \in R^n \times R^n$ be an arbitrary vector. Suppose that at each time moment $t \in [0, T]$ at the point $(x, \dot x) \in R^n \times R^n$ there exists such convex compact set $\underline \partial f(x,\dot x,t)$ $\subset R^n \times R^n$ that
\begin{equation}
\label{3} 
\frac{\partial f(x,\dot x,t)}{\partial g} = \lim_{\alpha \downarrow 0} \frac{1}{\alpha} \big(f(x+\alpha g_1, \dot x + \alpha g_2, t) - f(x,\dot x,t)\big) = \max_{v \in \underline \partial f(x,\dot x,t)} \langle v, g \rangle. \end{equation}

In this case the function $f(x,\dot x,t)$ is called subdifferentiable at the point $(x, \dot x)$ and the set $\underline \partial f(x,\dot x,t)$ is called the subdifferential of the function $f(x,\dot x,t)$ at the point $(x,\dot x)$. 

From expression (\ref{3}) one can see that at each $t \in [0, T]$ the following formula
\begin{equation}
\label{3''} f(x + \alpha g_1, \dot x + \alpha g_2, t) = f(x, \dot x, t) + \alpha  \frac{\partial f(x,\dot x,t)}{\partial g} + o(\alpha, x, \dot x, g, t), \end{equation}
$$\quad  \frac{o(\alpha, x, \dot x, g, t)}{\alpha} \rightarrow 0, \ \alpha \downarrow 0,$$
holds true.

If for each number $\varepsilon > 0$ there exist such numbers $\delta > 0$ and $\alpha_0 > 0$ that at $\overline g \in B_{\delta}(g)$ and $\alpha \in (0, \alpha_0)$ one has $| o(\alpha, x, \dot x, \overline g, t) | < \alpha \varepsilon $, then the function $f(x,\dot x,t)$ is called uniformly subdifferentiable at the point $(x, \dot x)$. Note \cite{DemyanovVasiliev} that if at each $t \in [0, T]$ the function $f(x,\dot x,t)$ is subdifferentiable at the point $(x, \dot x)$ and locally Lipschitz continuous in the vicinity of the point $(x, \dot x)$, then it is uniformly subdifferentiable at the point $(x, \dot x)$. If for the uniformly subdifferentiable function $f(x,\dot x,t)$ in expression (\ref{3''}) one has $\displaystyle{\frac{o(\alpha, x, \dot x, g, t)}{\alpha} \rightarrow 0}$, $\alpha \downarrow 0$, uniformly in $t \in [0, T]$, then such a function is called absolutely uniformly subdifferentiable. 

Consider the set $C_n[0, T] \times P_n[0, T]$ with the norm $L_2^n [0, T] \times L_2^n [0, T]$. Let $g = [g_1, g_2] \in C_n[0, T] \times P_n[0, T]$ be an arbitrary vector-function. Suppose that at the point $(x, z) \in C_n[0, T] \times P_n[0, T]$ there exists such a convex weakly* compact set $\underline \partial {I(x, z)} \subset \big( C_n[0, T] \times P_n[0, T], || \cdot ||_{L_2^n [0, T] \times L_2^n [0, T] } \big) ^*$ that  
\begin{equation}
\label{3'} 
\frac{\partial I(x, z)}{\partial g} = \lim_{\alpha \downarrow 0} \frac{1}{\alpha} \big(I(x+\alpha g_1, z+\alpha g_2) - I(x, z)\big) =  \max_{v \in \underline \partial I(x, z)} v(g). \end{equation}

In this case the functional $I(x,z)$ is called subdifferentiable at the point $(x, z)$, and the set $\underline \partial {I(x, z)}$ is called the subdifferential of the functional $I(x, z)$ at the point $(x,z)$. 

From expression (\ref{3'}) one can see that the following formula
\begin{equation}
\label{3'''}I(x + \alpha g_1, z + \alpha g_2) = I(x, z) + \alpha  \frac{\partial I(x,z)}{\partial g} + o(\alpha, x, z, g),\end{equation} $$ \quad  \frac{o(\alpha, x, z, g)}{\alpha} \rightarrow 0, \ \alpha \downarrow 0,$$
holds true.

So, it is required to find such vector-function $x^{*} \in C_n[0,T]$, which minimizes functional (\ref{1}) and satisfies boundary conditions (\ref{2}). Assume that there exists such a solution. The difference between this problem and the classical one of  variational calculus is that the integrand in the problem under consideration may not be smooth and be only subdifferentiable. 

Although the paper considers only continuous trajectories with a piecewise continuous and bounded derivative (this is due to the possibility of finding such trajectories in practice) and integrands with sufficiently ``good'' properties, let us give some known existence theorems for the considered problem with the solution in a class of all absolutely continuous functions and less burdensome constraints on the integrand.  In literature the scheme for proving the existence of a given variational problem solution consists of two main stages: proving compactness of a certain level set and proving lower semicontinuity of the considered functional in some topology. Let us give rather general results contained in papers \cite{Ioffe1}, \cite{Rock1}. The first one is formulated in terms of the integrand of the initial functional, while the second one uses the properties of the function conjugate to the integrand. 

 {\bf Theorem 1.1} 
 Let the integrand $f(x,y,t)$ in functional (\ref{1}) satisfy the following conditions:
 
 1) the function $f(x,y,t)$ maps the space $R^n \times R^n \times [0, T]$ to the interval $(-\infty, \infty]$;
 
 2) the function $f(x,y,t)$ is measurable with respect to the sigma algebra generated in the space $R^n \times R^n \times [0, T] $ by the direct product of Borel measurable subsets of the space $R^n \times R^n$ and Lebesgue measurable subsets of the segment $[0, T]$;
 
 3) the function $f(x,y,t)$ is lower semicontinuous in $(x,y)$ at each fixed $t \in [0, T]$; 
 
 4) the function $f(x, y, t)$ is convex in $y$ at each fixed $t \in [0, T]$, $x \in R^n$;
 
 5) at each $x, y \in R^n$ and at almost every $t \in [0, T]$ the inequality 
 $$f(x,y,t) \geq p(\|y\|_{R^n}) - q(\|x\|_{R^n}) + r(t)$$
 holds true where 
 
 5a) $p(\omega)$ is an nonnegative convex function defined on the interval $[0, \infty)$, and $p(0) = 0$;
 
 5b) $q(\omega)$ is an nonnegative continuous and nondecreasing function on the interval $[0, \infty)$; 
 
 5c) $\frac {p(\omega)}{\omega} \rightarrow \infty$ if $\omega \rightarrow \infty$;
 
 5d) $p\left( \frac {2 \omega}{T} \right) - q(\omega_0 + \omega) \rightarrow \infty$ if $\omega \rightarrow \infty$ where $\omega_0 = \max(\|x_0\|_{R^n}, \|x_T\|_{R^n})$;
 
 5e) the function $r(t)$ is summable on $[0,T]$. 
  
 Then, if only for one absolutely continuous function $x(t)$ satisfying conditions (\ref{2}) integral (\ref{1}) is finite, then problem (\ref{1}), (\ref{2}) has a solution in a class of absolutely continuous functions.
 
 {\bf Theorem 1.2} 
 Let the integrand $f(x,y,t)$ in functional (\ref{1}) satisfy conditions 1), 2), 3), 4) of Theorem 1.1.
 
 Let also the function $h(x,w,t)$ conjugate to the integrand and defined by the formula
 $$h(x,w,t) = \sup_y\{ \langle w, y \rangle - f(x,y,t) \}$$
 satisfy the following growth condition
 $$h(x,w,t) \leq \mu(w, t) + \|x\|_{R^n}(\sigma(t) + \rho(t) \|w\|_{R^n})$$
 where functions $\sigma(t)$, $\rho(t)$ are finite, nonnegative and summable on $[0, T]$, the function $\mu(w, t)$ is finite and summable in $t \in [0, T]$ at each fixed $w \in R^n$.
 
Then if only for one absolutely continuous function $x(t)$ integral (\ref{1}) is finite, then problem (\ref{1}) has a solution in a class of absolutely continuous functions.

These theorems are presented here in a somewhat reduced form than in papers \cite{Ioffe1}, \cite{Rock1} in order to avoid introducing more general spaces and corresponding metrics than those considered in this paper. As can be seen from the formulations of these theorems, they use only the general properties of the functions included in the formulation of the problem, such as continuity, measurability, convexity, etc. It is easy to verify the fulfillment of these conditions in a wide number of cases. 

%
%

\section{Reduction to an Unconstrained Minimization Problem}
Construct the functional, taking into account all the restrictions in the formulation of the problem. Let $z(t) = \dot x (t)$ (as we have assumed, $z \in P_n[0,T]$), then by virtue of restriction on the initial state (see the first equality in formula (\ref{2}) from the section Statement of the Problem) we have \begin{equation} 
\label{6} \displaystyle{x(t) = x_0 + \int_0^t z(\tau) d \tau}. \end{equation}
 
Construct the following functional on the space $P_n[0, T]$
\begin{equation} \label {3.1} \overline I(z) = \overline J\Big(x_0 + \int_0^t z(\tau) d \tau\Big) + \lambda \psi(z) = \end{equation} 
$$= \int_0^T f\Big(x_0 + \int_0^t z(\tau) d \tau, z(t), t\Big) dt + \lambda \frac{1}{2} \left( x_0 + \int_0^T z(t) dt - x_T \right)^2. $$
In functional $\overline I(z)$ the penalty summand with some positive factor $\lambda$ takes into account restriction on the final state of the system (see the second equality in formula (\ref{2}) from the section Statement of the Problem). 

Transition to the ``space of derivatives'' $z \in P_n[0, T]$ has been used in many works of V. F. Demyanov and his students to study various variational and control problems. Under some natural additional assumptions (namely: the function $f(x, z, t)$ is continuous in $(x,z,t)$, absolutely uniformly subdifferentiable and the mapping $t \rightarrow \underline \partial f(x(t), z(t),  t)$ is upper semicontinuous) one can prove the subdifferentiability of the functional $\overline I (z)$ in the space $P_n[0, T]$ as a normed space with the norm $L_2^n [0, T] $. However, the subdifferential of this functional has a rather complicated structure, which makes it practically unsuitable for constructing numerical methods. Therefore, it is proposed to consider some modification of this functional, ``forcibly'' considering the points $z$ and $x$ to be ``independent'' variables. Since, in fact, there is relationship (\ref{6}) between these variables (which naturally means that the function $z (t)$ is a derivative of the function $x (t)$), let us take this into account by adding the corresponding (last) term when constructing the new functional on the space $C_n[0, T] \times P_n[0, T]$   
\begin{equation} 
\label{7}
I(x, z) =  J(x, z) + \lambda \psi(z) + \lambda \varphi(x, z) =  
\end{equation}
$$
= \int_0^T f (x(t), z(t), t) dt +$$ $$+ \lambda \frac{1}{2} \left( x_0 + \int_0^T z(t) dt - x_T \right)^2 + \lambda \frac{1}{2} \int_0^T \Big( x(t) - x_0 - \int_0^t z(\tau) d \tau \Big)^2 dt.
$$

Despite the fact that the dimension of functional $I(x, z)$ arguments is $n$ more the dimension of functional $\overline I(z)$ arguments, the structure of its subdifferential (in the space $C_n[0, T] \times P_n[0, T]$ as a normed space with the norm $L_2^n [0, T] \times L_2^n [0, T]$), as will be seen from what follows, is much simpler than the structure of the functional $\overline I(z)$ subdifferential. This will allow us to construct a numerical method for solving the original problem. 

It is known \cite{Vasil'ev} that when the value $\lambda$ is sufficiently large, the solution of problem (\ref{1}), (\ref{2}) is arbitrarily close (with regard to the metric $L_2^n[0, T]$) to the trajectory $\overline x (t) $
where $(\overline x, \overline z)$ is a point of the global minimum of functional (\ref{7}) with the fixed value $\overline \lambda$. So, finding an approximate solution of the original problem is reduced to minimizing functional (\ref{7}) on the space \linebreak $C_n [0, T] \times P_n [0, T]$. In practice one solves this problem for the fixed number~$\overline{\lambda}$. If the solution of this problem (at $\lambda$ = $\overline{\lambda}$) satisfies the constraints in the form of differential relation (\ref{6}) and right endpoint condition from (\ref{2}) with the required accuracy (i. e. the value of the functional $\psi + \varphi$ on this solution is sufficiently small), then the process terminates; otherwise, increase the value $\lambda$ and repeat the process with this new value.  

Thus, the initial problem has been reduced to finding the unconditional global minimum point of the functional $I(x, z)$ (for sufficiently large value $\overline{\lambda}$) on the space 
\begin{equation} 
\label{07}
X = \big( C_n[0, T] \times P_n[0, T], || \cdot ||_{L_2^n [0, T] \times L_2^n [0, T]} \big).
\end{equation}

{\bf Remark 1.}
Note the following fact. Since, as is known, the space \linebreak $\big( C_n[0, T], || \cdot ||_{L_2^n [0, T]} \big)$ is everywhere dense in the space ${L_2^n [0, T]}$ and the space $\big( P_n[0, T], || \cdot ||_{L_2^n [0, T]} \big)$ is also everywhere dense in the space ${L_2^n [0, T]}$, then the space $X^*$ conjugate to the space $X$ (see (\ref{07})) is isometrically isomorphic to the space ${L_2^n [0, T] \times L_2^n [0, T]}$ \cite{KolmFom}; therefore, henceforth, we will identify these spaces ($X^*$ and ${L_2^n [0, T] \times L_2^n [0, T]}$).

{\bf Remark 2.}
If we consider the problem with the free right endpoint, then in formula (\ref{7}) one should put $\psi(z) = 0$ identically. If the minimized functional does not depend on the derivative $\dot x$, then in formula (\ref{7}) one should also put $\varphi(x,z) = 0$ identically. Since, by assumption, there exists a solution of the original problem, then the set of the trajectories $x$ satisfying the necessary minimum condition of the functional $I(x, z)$ (see the next section) will include those that satisfy boundary conditions  (\ref{2}). However, despite the fact that the minimized functional does not depend on the derivative in the case under consideration, one can also solve the problem ``fully'' by finding both unknowns in the pair $(x, \dot x)$ simultaneously, i.~e.  minimize functional (\ref{7}) including all its summands. Then (since, by assumption, there exists a solution of the original problem) among the trajectories $x$ satisfying the necessary minimum condition of the functional $I(x, z)$ there will be those that satisfy boundary conditions~(\ref{2}).  

\section{Minimum Conditions of the Functional $I(x, z)$}

In order to obtain the constructive minimum condition useful for constructing numerical methods for solving the posed problem, first, let us investigate the differential properties of the functional $I(x, z)$.

Using classical variation, it is easy to show Gateaux differentiability of the functional $\psi(z)$, we have 
$$
 \nabla \psi(z) = x_0 + \int_0^T z(t) dt - x_T.  
 $$
 
 Using classical variation and integrating by parts, it is also not difficult to check Gateaux differentiability of the functional $\varphi(x, z)$, we obtain
 $$\nabla \varphi(x, z, t) = \begin{pmatrix}
\displaystyle{ x(t) - x_0 - \int_0^t z(\tau) d\tau } \\
\displaystyle{ -\int_t^T \Big( x(\tau) - x_0 - \int_0^\tau z(s) ds \Big) d\tau  }
\end{pmatrix}. 
$$   
 
 
 Let us now study the differential properties of the functional \linebreak $\displaystyle{\int_0^T f(x(t), z(t), t) dt }$. Insofar as in this functional $x$ and $z$ are considered as independent variables, put $\xi(t) = (x(t), z(t))$ for brevity and prove the following theorem retaining the previous notation for the functional $J(x,z)$. 
 
 {\bf Theorem 2.}  
Consider the functional
 $$
 J(\xi) = \int_0^T f(\xi(t), t) dt, 
 $$
where $\xi \in C_n[0, T] \times P_n[0,T]$, the function $f(\xi, t)$ is continuous in $(\xi, t)$ and is absolutely uniformly subdifferentiable and its subdifferential is $\underline \partial f(\xi, t)$. Suppose also that the mapping $t \rightarrow \underline \partial f(\xi(t), t)$ is upper semicontinuous.  

Then the functional $J(\xi)$ is subdifferentiable, i. e. 
\begin{equation}
\label{10} 
\frac{\partial J(\xi)}{\partial g} = \lim_{\alpha \downarrow 0} \frac{1}{\alpha} \big(J(\xi+\alpha g) - J(\xi)\big) =  \max_{v \in \underline \partial J(\xi)} \int_0^T \langle v(t), g(t) \rangle dt,  \end{equation}

 where $g \in C_n[0, T] \times P_n[0,T]$ and the set $\underline \partial J(\xi)$ is defined as follows:
\begin{equation}
\label{12}  
 \underline \partial J(\xi) = \Big\{ v(t) \in L^{2n}_\infty[0, T] \ \big| \ v(t) \in  \underline \partial f(\xi(t), t) \ \forall t \in [0, T] \Big\}.
 \end{equation} 
  {\bf Proof.} 
In accordance with definition (\ref{3'}) of a subdifferentiable functional, to prove the theorem one has to check that:

1) the derivative of the functional $J(\xi)$ in the direction $g$ exists and is actually of form (\ref{10}), 

2) herewith, the set $\underline \partial J(\xi)$ is a convex and weakly* compact subset of the space $\big( C_n[0, T] \times P_n[0, T], || \cdot ||_{L_2^n [0, T] \times L_2^n [0, T]} \big)^*$. 

Prove statement 1). As the function $f(\xi, t)$ is subdifferentiable by assumption, then for every $g \in C_n[0, T] \times P_n[0, T]$ and for every $\alpha > 0$ we have
\begin{equation}
\label{16} 
J(\xi+\alpha g) - J(\xi) =  \int_0^T \max_{v \in \underline \partial f(\xi,t)} \langle v, \alpha g \rangle dt + \int_0^T {o(\alpha, \xi, g, t)} dt,
 \end{equation} 
  $$
 \frac{o(\alpha, \xi, g, t)}{\alpha} \rightarrow 0, \ \alpha \downarrow 0. 
 $$

As $\xi, g \in C_n[0, T] \times P_n[0, T]$ and the function $f(\xi, t)$ is continuous, one has that for each $\alpha > 0$ the functions $t \rightarrow f(\xi(t), t)$ and $t \rightarrow f(\xi(t) + \alpha g(t), t)$ belong to the space $L_\infty^1 [0,T]$.

Under the assumption made, the mapping $t \rightarrow \underline \partial f(\xi(t), t)$ is upper semicontinuous. Then due to the piecewise continuity of the function $g(t)$ and due to the continuity of the scalar product in its variables we obtain that  for each $\alpha >0$ the mapping $t \rightarrow \max_{v \in \underline \partial f(\xi(t),t)} \langle v, \alpha g(t) \rangle$ is upper semicontinuous~\cite{AubinFrankowska} and then it is also measurable \cite{FilippovBlagodatskih}. During the proof of statement 2) it will be shown that under the assumptions made, the set $\underline \partial f(\xi,t)$ is uniformly in $t \in [0, T]$ bounded; from here, taking into account the piecewise continuity of the function $g(t)$ it is easy to check that for each $\alpha >0$ the mapping $t \rightarrow \max_{v \in \underline \partial f(\xi(t),t)} \langle v, \alpha g(t) \rangle$ is also uniformly in $t \in [0, T]$ bounded. So we finally have that for each $\alpha >0$ the mapping $t \rightarrow \max_{v \in \underline \partial f(\xi(t),t)} \langle v, \alpha g(t) \rangle$ belongs to the space $L_\infty^1 [0,T]$. 

Then for every $\alpha > 0$ one has $t \rightarrow {o(\alpha, \xi(t), g(t), t)} \in  L_\infty^1 [0,T]$ and due to the absolutely uniformly subdifferentiability of the function $f(\xi, t)$ we have 
\begin{equation}
\label{11}  \frac{o(\alpha, \xi(t), g(t), t)}{\alpha} =: \frac{o(\alpha)}{\alpha} \rightarrow 0, \ \alpha \downarrow 0. \end{equation}

Consider the functional $\displaystyle{\int_0^T \max_{v \in \underline \partial f(\xi,t)} \langle v, \alpha g \rangle dt}$ in detail. For each $\alpha > 0$ and for each $t \in [0, T]$ we have the obvious equality
$$
\max_{v \in \underline \partial f(\xi,t)} \langle v(t), \alpha g(t) \rangle \geqslant \langle v(t), \alpha g(t) \rangle,
$$
where $v(t)$ is a measurable selector of the mapping $t \rightarrow \underline \partial f(\xi(t),t)$ (due to the noted boundedness property of the set $\underline\partial f(\xi, t)$ uniformly in $t \in [0, T]$ we have $v \in L_\infty^{2n} [0,T]$) and by virtue of formula (\ref{12}) for every $\alpha > 0$ one has the inequality
$$
\int_0^T \max_{v \in \underline \partial f(\xi,t)} \langle v, \alpha g \rangle dt \geqslant \max_{v \in \underline \partial J(\xi)} \int_0^T \langle v(t), \alpha g(t) \rangle dt. 
$$
As for every $\alpha > 0$ and for each $t \in [0, T]$ one has 
$$
\max_{v \in \underline \partial f(\xi,t)} \langle v(t), \alpha g(t) \rangle \in \Big\{ \langle v(t), \alpha g(t) \rangle \ \big| \ v(t) \in \underline \partial f(\xi(t),t) \Big\},
$$
and the set $\underline \partial f(\xi,t)$ is closed and bounded at each fixed $t$ by the definition of subdifferential and the mapping $t \rightarrow \underline \partial f(\xi(t),t)$ is upper semicontinuous by assumption and also because the scalar product is continuous in its arguments and the function $g(t)$ is piecewise continuous, then due to Filippov lemma \cite{Filippov} there exists such measurable selector $\overline{v}(t)$ of the mapping $t \rightarrow \underline \partial f(\xi(t),t)$ that for each $\alpha > 0$ and for each $t \in [0, T]$ we have 
$$
\max_{v \in \underline \partial f(\xi,t)} \langle v(t), \alpha g(t) \rangle = \langle \overline{v}(t), \alpha g(t) \rangle,
$$
so we have found the element $\overline{v}$ from the set $\underline \partial J(\xi)$ which brings the equality in the previous inequality. Thus, finally we obtain
\begin{equation}
\label{14} 
\int_0^T \max_{v \in \underline \partial f(\xi,t)} \langle v, \alpha g \rangle dt = \max_{v \in \underline \partial J(\xi)} \int_0^T \langle v(t), \alpha g(t) \rangle dt. 
\end{equation}


From (\ref{16}), (\ref{11}), (\ref{14}) we obtain expression (\ref{10}). 

Let us prove statement 2). The convexity of the set $\underline \partial J(\xi)$ immediately follows from the convexity of the set $\underline \partial f(\xi,t)$ at each fixed $t \in [0, T]$.

Prove the boundedness of the set $\underline \partial f(\xi,t)$ uniformly $t \in [0, T]$. Due to upper semicontinuity of the mapping $t \rightarrow \underline \partial f(\xi(t),t)$ at each $t \in [0, T]$ there exists such number $\delta(t)$ that under the condition $|\overline{t} - t| < \delta(t)$ the inclusion $\underline \partial f(\xi(\overline{t}),\overline{t}) \subset B_r(\underline \partial f(\xi({t}),{t}))$ holds true at $\overline{t} \in [0, T]$ where $r$ is some fixed finite positive number. The intervals $D_{\delta(t)}(t)$, $t \in [0, T]$, form open cover of the segment $[0, T]$, so by Heine-Borel lemma one can take a finite subcover from this cover. Hence, there exists such number $\delta > 0$ that for every $t \in [0, T]$ the inclusion $\underline \partial f(\xi(\overline{t}),\overline{t}) \subset B_r(\underline \partial f(\xi({t}),{t}))$ holds true once $|\overline{t} - t| < \delta$ and $\overline{t} \in [0, T]$. This means that for the segment $[0, T]$ there exists a finite partition $t_1 = 0, t_2, \dots, t_{N-1}, t_N = T$ with the diameter $\delta$ such that $\underline \partial f(\xi,{t}) \subset \bigcup\limits_{i=1}^N B_r(\underline \partial f(\xi({t_i}),{t_i}))$ for all $t \in [0, T]$. It remains to notice that the set $\bigcup\limits_{i=1}^N B_r(\underline \partial f(\xi({t_i}),{t_i}))$ is bounded due to the compactness of the set $\underline \partial f(\xi,{t})$ at each fixed $t \in [0, T]$.

Prove that the set $\underline \partial J(\xi)$ is weakly closed. As shown in statement~1) proof and at the beginning of statement 2) proof, the set $\underline \partial J(\xi)$ is convex and its elements $v$ belong to the space $L_\infty^{2n} [0, T]$. Then all the more the set $\underline \partial J(\xi)$ is a convex subset of the space $L_2^{2n} [0, T]$. Let us prove that the set $\underline \partial J(\xi)$ is closed in the weak topology of the space $L_2^{2n} [0, T]$. Let $\{v_n\}_{n=1}^{\infty}$ be the sequence of functions from the set $\underline \partial J(\xi)$ converging to the function $v^*$ in the strong topology of the space $L_2^{2n} [0, T]$. It is known \cite{Munroe} that this sequence has the subsequence $\{v_{n_k}\}_{n_k=1}^{\infty}$ converging pointwise to $v^{*}$ almost everywhere on $[0, T]$, i. e. there exists such subset $T' \subset [0, T]$ having the measure $T$ that for every point $t \in T'$ we have $v_{n_k}(t) \in \underline \partial f(\xi(t), t)$ and $v_{n_k}(t)$ converges to $v^{*}(t)$, $n_k = 1, 2, \dots$. But the set $\underline \partial f(\xi(t), t)$ is closed at each $t \in [0, T]$ by the definition of the subdifferential, hence for every $t \in T'$ we have $v^{*}(t) \in \underline \partial f(\xi(t), t)$. So the set $\underline \partial J(\xi)$ is closed in the strong topology of the space $L_2^{2n} [0, T]$, but it is also convex, so it is also closed in the weak topology of the space $L_2^{2n} [0, T]$ \cite{DunfordSchwartz}.

Recall that by virtue of Remark 1 it is sufficient to consider the space $L_2^{2n}[0, T]$. The weak* compactness of the set $\underline \partial J(\xi)$ in the space $L_2^{2n}[0, T]$ follows from its weak compactness (in $L_2^{2n}[0, T]$) by virtue of these topologies definitions \cite{KolmFom}. The space $L_2^{2n}[0, T]$ is reflexive \cite{DunfordSchwartz}, so the set there is weakly compact if and only if it is bounded in norm and weakly closed \cite{DunfordSchwartz} in this space. These required properties have been proved in the previous two paragraphs. The theorem is proved. 

Thus, as one can see from Theorem 2, the subdifferential of the functional~$J(\xi)$ is completely defined by the subdifferential of its integrand (at each time moment $t \in [0, T]$). So in order to calculate the subdifferential of the functional~$J(x, z)$, one has to  calculate the set $\underline \partial f(x, \dot x, t)$ for each $t \in [0,T]$ via subdifferential calculus \cite{DemyanovVasiliev} . Book \cite{DemyanovVasiliev} contains a detailed description of the rules for calculating the subdifferential for a rich class of functions. Let us recall some of these rules which are required while calculating the subdifferential of the functional $I(x, z)$.

 Let $\xi \in R^l$. If the function $\varphi(\xi)$ is subdifferentiable at the point $\xi_0 \in R^l$ and $\lambda$ is some nonnegative number, then one has
$$ \underline\partial ( \lambda  \varphi(\xi_0) ) = \lambda \, \underline\partial \varphi(\xi_0).$$
If the function $\varphi(\xi)$ is differentiable at the point $\xi_0 \in R^l$, then its subdifferential at this point is expressed by the formula $$\underline\partial \varphi(\xi_0)  = \varphi'(\xi_0), $$ where $\varphi'(\xi_0)$ is a gradient of the function $\varphi(\xi)$ at the point $\xi_0$.
Note that the subdifferential of the sum of the finite number of subdifferentiable functions is the that of the subdifferentials of the summands, i. e. if the functions $\varphi_k(\xi)$, $k = \overline{1, r}$, are subdifferentiable at the point $\xi_0 \in R^l$, then the subdifferential of the function $\varphi(\xi) = \sum_{k=1}^r \varphi_k(\xi)$ at this point is expressed by the formula $$ \underline\partial \varphi(\xi_0) =  \sum_{k=1}^r \underline\partial \varphi_k(\xi_0).$$
The subdifferential of the maximum of the finite number of continuously differentiable functions is the convex hull of the active functions gradients, i. e. if $\varphi(\xi_0) = \max \{\varphi_1(\xi_0), \dots, \varphi_r(\xi_0)\}$, then one has
 $$\underline\partial \varphi(\xi_0) = \mathrm{co}\{\varphi'_k(\xi_0)\}, \quad \text{if} \ \varphi_k(\xi_0) = \varphi(\xi_0), \ k \in \{1..r\},$$
 where $\varphi_k'(\xi_0)$ is a gradient of the function $\varphi_k(\xi)$ at the point $\xi_0$, $k \in \{1..r\}$.


Using formula (\ref{12}) and these rules of subdifferential calculus, one obtains the final expression for calculating the subdifferential of the functional $I(x, z)$ at the point $(x, z)$
\begin{equation}
\label{16'}
\underline \partial I(x, z) = \sum_{k=1}^3 \underline \partial I_k(x, z),
\end{equation} 
where formally $I_1(x, z) = J(x, z)$, $I_2(x, z) = \lambda \psi(z)$, $I_3(x, z) = \lambda \varphi(x, z)$.

The known necessary minimum condition of the functional $I(x, z)$ at the point $(\overline x, \overline z)$ in terms of subdifferential is as follows \cite{Dolgopolik1} 
$$
0_{2n} \in \underline \partial I (\overline x, \overline z),
$$
where $0_{2n}$ is a zero element of the space $L_2^{2n} [0, T]$.
Hence, we conclude that the following theorem is true. 

{\bf Theorem 3.} 
For the point $(\overline x, \overline z)$ to minimize functional (\ref{7}), it is necessary that for almost every $t \in [0, T]$ the inclusion
\begin{equation}
\label{17}
\bold 0_{2n} \in \underline \partial I (\overline x(t), \overline z(t))
\end{equation}
is satisfied, where $\bold 0_{2n}$ is a zero element of the space $R^{2n}$, and the expression for the subdifferential $\underline \partial I (x, z)$ is given by formula (\ref{16'}). 

{\bf Remark 3.} 
Theorem 3 contains a constructive minimum condition since on its basis it is possible to construct the subdifferential descent method which is described in the next section. Although the principle algorithm of this method is well-known, its application to functional $\overline I(z)$ in formula (\ref{3.1}) is impossible in practice since the structure of this functional subdifferential is too complicated and it is unclear how to solve the important subproblems of the algorithm. So the key idea of this paper is ``forcibly'' considering the points $z$ and $x$ to be ``independent'' variables and to applicate the subdifferential descent method to the constructed functional $I(x,z)$ in formula (\ref{7}). (See additional explanations after formulas (\ref{3.1}) and (\ref{7}) in the section Reduction to an Unconstrained Minimization Problem.) This idea makes it possible (see the next section) to solve each of the arising subproblems of the method via known effective algorithms. 

\section {The Subdifferential Descent Method}

Describe the following subdifferential descent algorithm for finding stationary points of the functional $I(x, z)$. 

Fix the arbitrary initial point $(x_{(1)}, z_{(1)}) \in C_n[0, T] \times P_n[0, T]$. Let the point $(x_{(k)}, z_{(k)}) \in C_n[0, T] \times P_n[0, T]$ be already constructed. If minimum condition (\ref{17}) is satisfied (in practice with some fixed accuracy $\overline{\varepsilon}$), then the point $(x_{(k)}, z_{(k)})$ is a stationary point of the functional $I(x, z)$ and the process terminates. Otherwise, put
$$
(x_{(k+1)}, z_{(k+1)}) = (x_{(k)}, z_{(k)}) + \gamma_{(k)} G\big(x_{(k)}, z_{(k)}\big),
$$
where the vector-function $G\big(x_{(k)}, z_{(k)}\big)$ is a subdifferential descent direction of the functional $I(x, z)$ at the point $(x_{(k)}, z_{(k)})$, and the value $\gamma_{(k)}$ is a solution of the following one-dimensional problem
\begin{equation}
\label{18}
\min_{\gamma \geqslant 0} I \Big( (x_{(k)}, z_{(k)}) + \gamma G\big(x_{(k)}, z_{(k)}\big) \Big) = I \Big( (x_{(k)}, z_{(k)}) + \gamma_{(k)} G\big(x_{(k)}, z_{(k)}\big) \Big). 
\end{equation}
Then, as it will be shown in this section, \begin{equation}
\label{18'}
I \big(x_{(k+1)}, z_{(k+1)}\big) < I \big(x_{(k)}, z_{(k)}\big). \end{equation}

As seen from this algorithm, one has to solve three subproblems in order to realize the $k$-th iteration. The first problem is calculating the subdifferential of the functional $I(x, z)$ at the point $(x_{(k)}, z_{(k)})$. With the help of subdifferential calculus rules the solution of this problem is obtained in formula (\ref{16'}). The second problem is finding the subdifferential descent direction $G\big(x_{(k)}, z_{(k)}\big)$; two next paragraphs are devoted to this problem. Finally, the third problem is one-dimensional minimization (\ref{18}); there exist many effective methods \cite{Vasil'ev} to solve this problem. 

In order to find the vector-function $G\big(x_{(k)}, z_{(k)}\big)$, consider the problem
\begin{equation}
\label{19}
\min_{v \in \underline{\partial} I(x_{(k)}, z_{(k)})} ||v||^2_{L_2^{n}[0,T] \times L_2^{n}[0, T]}  = \min_{v \in \underline{\partial} I(x_{(k)}, z_{(k)}) } \int_0^T v^2(t)  dt. 
\end{equation}
Denote $\overline{v}_{(k)}$ the solution of this problem. The vector-function $\overline{v}_{(k)}(t)$, of course, depends on the point $(x_{(k)}, z_{(k)})$, but we omit this dependence in the notation for brevity. Then the vector-function $$G\big(x_{(k)}(t), z_{(k)}(t), t\big) = -\frac{\overline{v}_{(k)}\big(x_{(k)}(t), z_{(k)}(t), t\big)}{||\overline{v}_{(k)}||_{L_2^{2n}[0,T] }}$$ is a subdifferential descent direction of the functional $I(x, z)$ at the point $(x_{(k)}, z_{(k)})$. Recall that we are seeking the direction $G\big(x_{(k)}, z_{(k)}\big)$  in the case when the point $(x_{(k)}, z_{(k)})$ does not satisfy minimum condition in Theorem 3, so $||\overline{v}_{(k)}||_{L_2^{2n}[0,T]} > 0$.    

Note that we have the equality
$$\frac{\partial I(x_k, z_k)}{\partial G(x_k, z_k)} = \max_{v \in \partial I(x_k, z_k)} \int_0^T \left \langle v(t), G(x_k(t), z_k(t), t) \right \rangle dt =$$ $$=\max_{v \in \partial I(x_k, z_k)} \int_0^T \left \langle v(t), \frac{-\overline{v}_{(k)}(t)}{||\overline{v}_{(k)}||_{L_2^{2n}[0, T]}} \right \rangle dt =$$ 
$$= \frac{-1}{||\overline{v}_{(k)}||_{L_2^{2n}[0, T]}} \left( -\max_{v \in \partial I(x_k, z_k)} \int_0^T \left \langle -v(t), {\overline{v}_{(k)}(t)} \right \rangle \right) dt =$$ 
$$ = \frac{-1}{||\overline{v}_{(k)}||_{L_2^{2n}[0, T]}} \left( \min_{v \in \partial I(x_k, z_k)} \int_0^T \left \langle v(t), {\overline{v}_{(k)}(t)} \right \rangle \right) dt = -||\overline{v}_{(k)}||_{L_2^{2n}[0,T]},$$
which considering (\ref{3'''}) and the inequality $||\overline{v}_{(k)}||_{L_2^{2n}[0,T]} > 0$ implies (\ref{18'}). 

It is easy to check that in this case the solution of this problem is such selector of the multivalued mapping $t \rightarrow \underline{\partial} I\big(x_{(k)}(t), z_{(k)}(t), t \big)$ which minimizes the distance from zero to the set $\underline{\partial} I\big(x_{(k)}(t), z_{(k)}(t), t )$ at each time moment $t \in [0, T]$. In other words, to solve problem (\ref{19}) means to solve the following problem  
\begin{equation}
\label{20}
\min_{v(t) \in \underline{\partial} I(x_{(k)}(t), z_{(k)}(t), t ) } v^2(t)
\end{equation}
for each $t \in [0, T]$.
Actually, for every $t \in [0, T]$ we have the obvious inequality
$$
\min_{v \in \underline{\partial} I(x_{(k)}(t), z_{(k)}(t), t ) }  v^2(t) \leqslant v^2(t) ,
$$
where $v(t)$ is a measurable selector of the mapping $t \rightarrow \underline{\partial} I\big(x_{(k)}(t), z_{(k)}(t), t \big)$ (by virtue of the noted property of the set $\underline{\partial} I\big(x_{(k)}(t), z_{(k)}(t), t \big)$ boundedness uniformly in $t \in [0, T]$ we have $v \in L_\infty^{2n} [0,T]$), then we obtain the inequality
$$
\int_0^T \min_{v \in \underline{\partial} I(x_{(k)}(t), z_{(k)}(t), t ) }  v^2(t) dt \leqslant \min_{v \in \underline{\partial} I(x_{(k)}, z_{(k)}) } \int_0^T v^2(t)  dt. 
$$
Insofar as for every $t \in [0, T]$ we have 
$$
\min_{v \in \underline{\partial} I(x_{(k)}(t), z_{(k)}(t), t ) }  v^2(t) \in \Big\{ v^2(t) \ \big| \ v(t) \in \underline{\partial} I(x_{(k)}(t), z_{(k)}(t), t ) \Big\}
$$
and the set $\underline{\partial} I\big(x_{(k)}(t), z_{(k)}(t), t \big) $ is closed and bounded at every fixed $t$ by definition of the subdifferential and the mapping $t \rightarrow \underline{\partial} I\big(x_{(k)}(t), z_{(k)}(t), t \big) $ is upper semicontinuous by assumption and besides, the norm is continuous in its argument, then due to Filippov lemma \cite{Filippov} there exists such a measurable selector $\overline{v}_k(t)$ of the mapping $t \rightarrow \underline{\partial} I\big(x_{(k)}(t), z_{(k)}(t), t \big) $ that for every $t \in [0, T]$ one obtains 
$$
\min_{v \in \underline{\partial} I(x_{(k)}(t), z_{(k)}(t), t ) }  v^2(t) = \overline{v}_k^2(t),
$$
so we have found the element $\overline{v}_k$ of the set $\underline{\partial} I\big(x_{(k)}, z_{(k)}\big) $ which brings the equality to the previous inequality. Hence, finally we obtain
$$
\int_0^T \min_{v \in \underline{\partial} I(x_{(k)}(t), z_{(k)}(t), t ) }  v^2(t) dt  = \min_{v \in \underline{\partial} I(x_{(k)}, z_{(k)}) } \int_0^T v^2(t)  dt. 
$$

Problem (\ref{20}) at each fixed $t \in [0, T]$ is a finite-dimensional problem of finding the distance from zero to a convex compact (the subdifferential). This problem can be effectively solved for a wide class of functions; the next paragraph describes its solution. In practice one makes a (uniform) partition of the interval $[0, T]$, and this problem is solved for every point of the partition, i. e. one has to calculate $G\big(x_{(k)}(t_i), z_{(k)}(t_i), t_i \big)$, where $t_i \in [0, T]$, $i = \overline{1, N}$, are the points of discretization (see notation in Lemma 1 below). Under some natural additional assumption Lemma 1 below guarantees that the vector-function obtained with the help of piecewise linear interpolation of the subdifferential descent directions evaluated at every point of such partition of the interval $[0, T]$ converges to the sought vector-function $G\big(x_{(k)}(t), z_{(k)}(t), t\big)$ in the space $L_2^{2n}[0,T]$ when the discretization rank tends to infinity.    

As noted in the previous paragraph, for the algorithm realization it is required to find the distance from zero to the subdifferntial of the functional $I (x(t), z(t))$ at each moment of time of a (uniform) partition of the interval $[0, T]$. Let us discuss some methods for solving this subproblem (for a wide class of functions) for the fixed time moment $t \in [0, T]$. It is known \cite{DemyanovVasiliev} that in many practical cases the subdifferential $\underline \partial I (x(t), z(t))$ is the convex polyhedron $A(t) \subset R^{2n}$. For example, if the integrand is a maximum of the finite number of continuously differentiable functions, then the subdifferential $\underline \partial I (x(t), z(t))$ is a convex polyhedron at each  $t \in [0,T]$. Herewith, of course, the set $A(t)$ depends on the point $(x, z)$. We will omit this dependence in the notation in this paragraph for simplicity. This problem of finding the Euclidean distance from a point to a convex polyhedron can be effectively solved by various methods (see, e. g., \cite{DemyanovMal}, \cite{wolfe}). In a more general case the subdifferential at each moment $t \in [0, T]$ of time may be a convex compact set (for example, if the integrand depends on the norm of some coordinates of the vector-functions $x(t)$, $z(t)$, then the subdifferential at some points $t \in [0, T]$ may be an ellipsoid (with its interior points), lying in some subspace of the space $R^{2n}$). In this case it is required to solve the problem of finding the Euclidean distance from a point to a convex compact set, and if (for example) ellipsoids are considered, then some methods for solving this problem can be found in  \cite{dolgtam}. 

Prove one lemma with a simple condition, which on the one hand, is rather natural for applications and on the other hand, guarantees that the function $L(t)$ obtained with the help of piecewise linear interpolation of the sought function $p \in L^1_{\infty}[0, T]$ converges to this function in the space $L^1_2[0, T]$ when the rank of a (uniform) partition of the interval $[0, T]$ tends to infinity.   

 {\bf Lemma 1.}
 Let the function $p \in L^1_{\infty}[0, T]$ satisfy the following condition: for every $\overline{\delta} > 0$ the function $p(t)$ is piecewise continuous on the set $[0, T]$ with the exception of only the finite number of the intervals $\big(\overline t_1(\overline{\delta}), \overline t_2(\overline{\delta})\big),$ $\dots$, $\big(\overline t_{r}(\overline{\delta}), \overline t_{r+1}(\overline{\delta})\big)$ whose union length does not exceed the number $\overline{\delta}$.  
 
 Choose the (uniform) finite splitting $t_1 = 0, t_2, \dots, t_{N-1}, t_N = T$ of the interval $[0, T]$ and calculate the values $p(t_i)$, $i = \overline{1, N}$, at these points. Let $L(t)$ be the function obtained with the help of piecewise linear interpolation with the nodes $(t_i, p(t_i))$, $i = \overline{1, N}$. Then for every $\varepsilon > 0$ there exists such number $\overline{N}(\varepsilon)$ that for every $N > \overline{N}(\varepsilon)$ one has $||L - p||^2_{L^1_2[0,T]} \leqslant \varepsilon$.  

  {\bf Proof.} Denote $M(\overline{\delta}) := \bigcup\limits_{k=1}^{r} \big(\overline t_{k}(\overline{\delta}), \overline t_{k+1}(\overline{\delta})\big)$. We have
  $$||L - p||^2_{L^1_2[0,T]} = \int_{M(\overline{\delta})} \big (L(t) - p(t) \big)^2 dt + \int_{[0, T] \setminus M(\overline{\delta})} \big (L(t) - p(t) \big)^2 dt.$$
   Fix the arbitrary number $\varepsilon > 0$. By lemma condition the function $p(t)$ is bounded, the function $L(t)$ is also bounded by construction for all (uniform) finite partitions of the interval $[0, T]$. Hence, there exists such $\overline\delta(\varepsilon)$ that the first summand does not exceed the value ${\varepsilon}/{2}$ for all (uniform) finite partitions of the interval $[0, T]$. As assumed, the function $p(t)$ is piecewise continuous and bounded on the set $[0, T] \setminus M(\overline{\delta}(\varepsilon))$, then there exists \cite{Ryab} such number $\overline{N}(\varepsilon)$ that for every (uniform) finite partition of the interval $[0, T]$ of the rank $N > \overline{N}(\varepsilon)$ the second summand (with such $\overline{\delta}(\varepsilon))$ does not exceed the value ${\varepsilon}/{2}$. This implies the proof of the lemma. 
   
   {\bf Remark 4.}   
   The problem of a rigorous proof of the above method convergence is rather complicated and remains open; it is beyond the scope of this paper. The convergence of some modifications (related to the choice of a descent step and a descent direction from the set of subgradients) of the subdifferential descent method described in this section was studied in the finite-dimensional case in papers \cite{DemyanovVasiliev}, \cite{DemyanovMal}. Strictly speaking, in presented paper only the problem of finding the direction of the steepest (subdifferential) descent in the problem posed is completely solved. The examples below show the adequacy of the method used; nevertheless, as has been just noted, its convergence (in whatever sense) requires additional rigorous justification. 

\section{Numerical Examples}
   
   In this section the examples of the subdifferential descent method implementation are presented. These particular examples are chosen in order to demonstrate the described method processing in some standard cases considering such typical subdifferential functions as modules, square roots, maxima of continuously differentiable functions, etc. The stopping criteria of the algorithm was the inequality  $||\overline v_{(k)}||^2_{L_2^{n}[0,T] \times L_2^{n}[0, T]} \leqslant \overline \varepsilon$ (see problem (\ref{19})). In different examples the value $ \overline \varepsilon $ was taken equal to $3 \times 10^{-2}$ --- $9 \times 10^{-2}$. Such a choice of accuracy is due to a compromise between the permissible for practice accuracy of the optimal value of the considered functional and a not very great number of iterations. Herewith, the error of the minimized functional and the restrictions on the right endpoint in the examples below did not exceed the value $10^{-3}$ --- $5 \times 10^{-3}$ (in those examples where it was possible to compare the values obtained with the known solution).   
    
\textbf{Example 1.}
Consider minimization of the simplest functional
$$\overline{J}(x) = \int_0^1 |x(t)| dt,$$ $$x(0) = 0,$$
with the only obvious solution $x^*(t) = 0$ $\forall t \in [0, 1]$ and $\overline J (x^*) = 0$. In accordance with Remark~2, the functionals $\psi(z)$ and $\varphi(x, z)$ are absent here. Take $x_{(1)}(t) = 2t - 1$ as the initial point and discretize the segment $[0, 1]$ with rank two (i. e. consider the points $0$, $0.5$, $1$ for further subdifferential descent direction interpolation). In accordance with the paper algorithm, separately calculate the descent directions at these points. Then at the point $t_ 1 = 0$ the function $|x_{(1)}(t)|$ is differentiable, hence its subdifferential is of the form $ \underline{\partial} x_{(1)}(0) = \{-1\} $, find the distance from zero to the set $\{-1\}$ and obtain the subdifferential descent direction $G(x_{(1)}, 0) = 1$. In a similar way we have $G(x_{(1)}, 1) = -1$. At the point $t_2 = 0.5$ the function $|x_{(1)}(t)|$ is subdifferentiable, hence its subdifferential is of the form $ \underline{\partial} x_{(1)}(0.5) = [-1, 1]$, find the distance from zero to the set $[-1, 1]$ and obtain the subdifferential descent direction $G(x_{(1)}, 0.5) = 0$. Interpolating with the nodes $(0, 1)$, $(0.5, 0)$, $(1, -1)$, obtain the subdifferential descent direction of the functional $\overline{J}$ at the point $x_{(1)}$, namely $G(x_{(1)}) = -2t+1$. Construct the next point $x_{(2)}(t) = 2t-1 + \gamma (-2t+1)$ and solving the one-dimensional problem $\min_{\gamma \geqslant 0} \displaystyle{\int_0^1 |x_{(2)}(t)| dt}$, we have $\gamma_{(1)} = 1$, hence $x_{(2)}(t) = 0 $ $\forall t \in [0, 1]$, i. e. in this case the method leads to the exact solution in one step. Of course, the initial point and the discretization rank are artificially chosen here in order to demonstrate the essence of the method. If we take a different initial point and some other discretization rank, then the solution will not be obtained (in general case) in a finite number of steps. 

\textbf{Example 2.}
Consider minimization of the functional 
$$\overline{J}(x) = \int_0^1 \left|x(t) - \max\left\{t-0.5, 0\right\} \right| dt, $$
$$x(0) = 0, $$
with the only obvious solution $x^*(t) = \max\{t-0.5, 0\}$ $\forall t \in [0, 1]$ and $\overline J (x^*) =~0$. In accordance with Remark 2, the functionals $\psi(z)$ and $\varphi(x, z)$ are absent here. Take $x_{(1)} = 2t - 1$ as the initial point, then $I(x_{(1)}) = 0.375$. As the iteration number increased, the discretization rank gradually increased during the solution of the auxiliary problem of finding the direction of the subdifferential descent described in the algorithm and in the end the discretization step was equal to $10^{-1}$. At the 28-th iteration the point $x_{(28)} = $ $$ 14.1565t^5-13.6885t^4+3.7699t^3-0.0789t^2-0.0739t+0.0049, \quad 0 \leqslant t < 0.5,$$ $$6.0666t^5-19.4749t^4+23.4983t^3-12.9012t^2+3.9828t-0.6695, \quad 0.5 \leqslant  t \leqslant 1,$$ was obtained and the value of the functional $\overline J (x_{(28)}) = I(x_{(28)}) \approx 0.00116$, so the error does not exceed the value $10^{-3}$. For the convenience of presentation, the Lagrange interpolation polynomial has been given, which quite accurately approximates (that is, the interpolation error does not affect the value of the functional presented with a given accuracy but (insignificantly) affects the given value of the norm of the smallest subgradient) the resulting trajectory. Herewith, $||\overline{v}_{(28)}||_{L_2^1}[0,T] \approx 0.032$. 

  \textbf{Example 3.}
Minimize the functional
$$\overline{J}(x) = \int_0^1 \max\left\{\dot{x}_1^2(t) - x_1^2(t) - 2 t x_1(t), x_2(t)\right\} dt, $$
$$x_1(0) = 0,  \quad x_2(0) = 0,$$
$$x_1(1) = 0,  \quad x_2(1) = 0.$$
So, one has to minimize the functional 
   $$I(x,z) = \int_0^1 \max\left\{\dot{x}_1^2(t) - x_1^2(t) - 2 t x_1(t), x_2(t)\right\} dt + $$
   $$ + \lambda \frac{1}{2} \left( \int_0^1 z_1(t) dt  \right)^2 + \lambda \frac{1}{2} \left( \int_0^1 z_2(t) dt  \right)^2 + \lambda \frac{1}{2} \int_0^1 \Big( x(t) -  \int_0^t z(\tau) d \tau \Big)^2 dt.$$

  The point $(x_{(1)}, z_{(1)}) = (0, 0, 0, 0, 0, 0)' $ was taken as the initial one, and $I(x_{(1)}, z_{(1)}) = 0$.  As the iteration number increased, the discretization rank gradually increased during the solution of the auxiliary problem of finding the direction of the subdifferential descent described in the algorithm and in the end the discretization step was equal to $5 \times 10^{-2}$; the penalty parameter value also increased and in the end we had $\lambda = 300$. At the 56-th iteration the point $(x_{(56)}, z_{(56)})$ was constructed with the following coordinates:
$$ x_1 = 27.83995 t^5 - 18.272210 t^4 + 4.163818 t^3 - 0.47695 t^2 + 0.153284 t, \quad 0 \leq t \leq 0.25, $$
$$ 0.763217 t^5 -0.923994 t^3 + 0.457991 t^2 + 0.018469 t + 0.009833, \quad 0.25 \leq t \leq 0.5, $$
$$ 2.02681 t^5 - 3.01454 t^4 - 0.62533 t^3 + 3.13597 t^2 - 1.81301 t + 0.36767, \quad 0.5 \leq t \leq 0.75, $$ $$-0.155737 t^4 + 0.420485 t^2 - 0.429341 t + 0.169993, \quad 0.75 \leq t \leq 1, $$ $$ x_2 = 1.934618 t^5 - 2.059840 t^4 + 0.821448 t^3 - 0.169958 t^2 - 0.073440 t, \quad 0 \leq t \leq 0.5, $$
$$0.463557 t^4 - 1.012055 t^3 + 1.092861 t^2 - 0.647782 t + 0.103397, \quad 0.5 \leq t \leq 1, $$
$$z_1 = -0.043907 t^4 + 0.217162 t^3 - 0.337550 t^2 - 0.126702 t + 0.132179, $$
$$z_2 = -0.035801 t^4+0.661239 t^3-0.132123 t^2 - 0.066889 t - 0.080925, $$ and the functional value $I(x_{(56)}, z_{(56)}) \approx -0.02175$, $x_1(1) \approx 0.0054$, $x_2(1)  \approx 0.$ For the convenience of presentation, the function has been given that consists of Lagrange interpolation polynomials separately calculated on several time intervals and approximating quite accurately (that is the interpolation error does not affect the value of the functional presented with a given accuracy but (insignificantly) affects the given value of the norm of the smallest subgradient) the resulting trajectory. Herewith, $||\overline{v}_{(56)}||_{L_2^{2}[0,T] \times L_2^{2}[0, T]} \approx 0.0915$. 

Insofar as in fact one has $z(t) = \dot x(t)$, $t \in [0,1]$, then the point $(x_{(56)}, z_{(56)})$ may be taken with the coordinates $$ x_1 = 27.83995 t^5 - 18.272210 t^4 + 4.163818 t^3 - 0.476947 t^2 + 0.153284 t, \quad 0 \leq t \leq 0.25, $$
$$ 0.763217 t^5 -0.923994 t^3 + 0.457991 t^2 + 0.018469 t + 0.009833, \quad 0.25 \leq t \leq 0.5, $$
$$ 2.02681 t^5 - 3.01454 t^4 - 0.625326 t^3 + 3.13597 t^2 - 1.81301 t + 0.36767, \quad 0.5 \leq t \leq 0.75, $$ $$-0.155737 t^4 + 0.420485 t^2 - 0.429341 t + 0.169993, \quad 0.75 \leq t \leq 1, $$ $$ x_2 = 1.934618 t^5 - 2.059840 t^4 + 0.821448 t^3 - 0.169958 t^2 - 0.073440 t, \quad 0 \leq t \leq 0.5, $$
$$0.463557 t^4 - 1.012055 t^3 + 1.092861 t^2 - 0.647782 t + 0.103397, \quad 0.5 \leq t \leq 1, $$
$$ z_1 = 139.199755 t^4-73.88839 t^3+12.491454 t^2-0.953894 t+0.153284, \quad 0 \leq t \leq 0.25, $$
$$ 3.816084 t^4-2.771981 t^2+0.915982 t+0.018469, \quad 0.25 < t \leq 0.5, $$
$$ 10.134052 t^4-12.058139 t^3-1.875979 t^2+6.271932 t-1.813011, \quad 0.5 < t \leq 0.75, $$ $$-0.622947 t^3+0.840969 t-0.429341, \quad 0.75 < t \leq 1, $$
$$ z_2 = 9.673089 t^4-8.239360 t^3+2.464344 t^2-0.339916 t-0.073440, \quad 0 \leq t \leq 0.5, $$
$$1.854229 t^3-3.036168 t^2+2.185722 t-0.647782, \quad 0.5 < t \leq 1, $$ then $ I(x_{(56)}, z_{(56)}) \approx -0.01827$ and $\overline{J}(x_{(56)}) \approx -0.02264$, $x_1(1) \approx 0.0054$, $x_2(1)  \approx 0 $, i. e. the error of the trajectory value at the right endpoint does not exceed the magnitude $5 \times 10^{-3}$.

In order to assess the adequacy of the result obtained, let us turn to the example considered in paper \cite{TamasyanDemyanov}. There the smooth functional \linebreak $\overline{J}_0(x_1) = \displaystyle{\int_0^1 \dot{x}_1^2(t) - x_1^2(t) - 2 t x_1(t) dt} $ has been minimized with the boundary conditions $x_1(0) = 0$, $x_1(1) = 0$ and the exact solution $\overline{x}_1(t) = \displaystyle{ \frac{\sin(t)}{\sin(1)} - t}$ has been given and the corresponding minimum value $\overline{J}_0(\overline{x}_1) = -0.02457$. It is, however, easy to observe that the minimum values of the functionals $\overline{J}_0(x_1)$ and $\overline{J}(x_1, x_2)$ (with the given boundary conditions) coincide. Indeed, since  $ \max\left\{\dot{x}_1^2(t) - x_1^2(t) - 2 t x_1(t), x_2(t)\right\} \geq \dot{x}_1^2(t) - x_1^2(t) - 2 t x_1(t)$ for arbitrary values of the variables in this inequality, then one obtains the inequality $\overline{J}(x_1, x_2) \geq \overline{J}_0(x_1) \geq \overline{J}_0(\overline{x}_1)$ $\forall x_1, x_2 \in C_1[0, T] $. Having taken the arbitrary trajectory $\overline{x}_2 \in C_1[0,1]$, $\overline{x}_2(0) = 0$, $\overline{x}_2(1) = 0$, such that \linebreak $\overline{x}_2(t) \leq \dot{\overline{x}}_1^2(t) - \overline{x}_1^2(t) - 2 t \overline{x}_1(t) $ (it is possible since the values of the function in the right-hand side of this inequality are positive at $t=0$ and at $t = 1$), we obtain the proof of the statement. Since the minimum value of the functional obtained with the help of the method of this paper is approximately equal to $-0.02264$, then one can see that the error of the functional value does not exceed the magnitude $2 \times 10^{-3}$. 

Despite the fact that, as just noted, the minimal value of the functional in this example may be achieved at any function $\overline x_2(t)$ satisfying the given boundary conditions and not exceeding (at every $t \in [0,1]$) the integrand of the functional $\overline{J}_0(x_1)$ at the most iterations of the method described in the paper the values of the function $x_2(t)$ and of the integrand of the functional~$\overline{J}_0(x_1)$ coincided on a subset of the interval $ [0, 1] $ of nonzero measure, i.~e. both of the functions under the maximum in the functional $\overline{J}(x)$ were active on this subset, hence the functional $\overline{J}(x)$ turned out to be nondifferentiable at these functions and during the operation of the paper method, the ``complete'' subdifferential of the given functional was calculated. Figure~1 with the functions \linebreak $z_1^2(t) - x_1^2(t)  -2 t x_1(t)$ and $x_2(t)$ depicted illustrates the described situation for some iterations. 
\newline

\text {\textbf{Table 1}\, Example 3}
$$
\begin{array}{c|c|c|c}
        k&I(x_{(k)}, z_{(k)})&||\overline v(x_{(k)}, z_{(k)})||&\lambda\\
  \hline & &  \\[-3mm]
     1 & 0 & 0.6733 & 20 \\
     20 & -0.01289 & 0.2274 & 100  \\
     40 & -0.01833 & 0.1389 & 200  \\
     56 &  -0.02175 & 0.0915 & 300 \\
\end{array}
$$

  \begin{figure*}[h!]
\begin{minipage}[h]{0.3\linewidth}
\center{\includegraphics[width=1\linewidth]{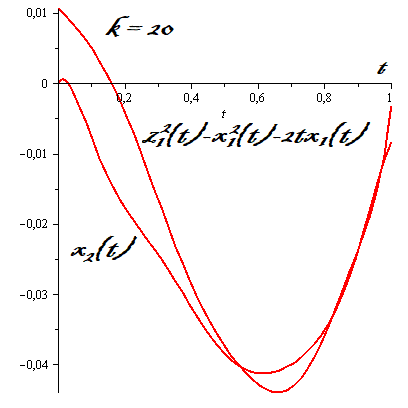} }
\end{minipage}
\hfill
\begin{minipage}[h]{0.3\linewidth}
\center{\includegraphics[width=1\linewidth]{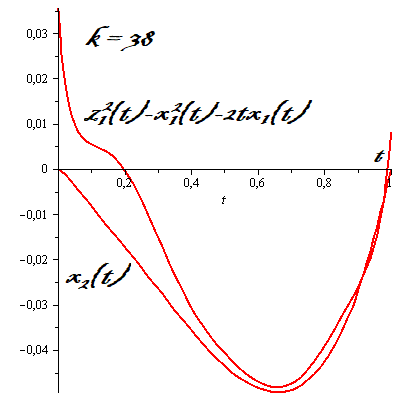} }
\end{minipage}
\hfill
\begin{minipage}[h]{0.3\linewidth}
\center{\includegraphics[width=1\linewidth]{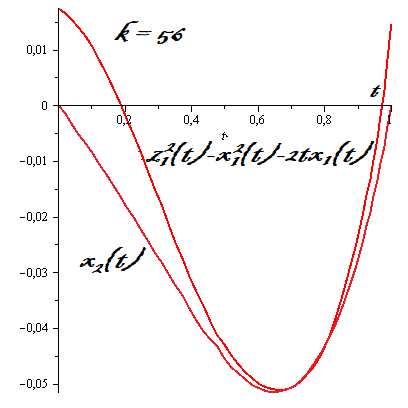} }
\end{minipage}
\caption{Example 3}
\label{ris:image2}
\end{figure*}

  \textbf{Example 4.}
Consider the minimization problem of the functional
$$\overline{J}(x) = \int_0^5 \sqrt{(\dot{x}_1(t) - 1)^2 + x_2^2(t)} + (x_1(t)-x_3(t)-\sin(t))^2 dt, $$
$$x_1(0) = 0,  \quad x_2(0) = 0, \quad x_3(0) = 0,$$
with the only obvious solution $x_1^*(t) = t$, $x_2^*(t)=0$, $x_3^*(t)=t-\sin(t)$, $t \in [0, 5]$, and $J(x^*) = 0$. In accordance with Remark 2, the functional $\psi(z)$ here is absent. So, it is required to minimize the functional
   $$I(x,z) = \int_0^5 \sqrt{(\dot{x}_1(t) - 1)^2 + x_2^2(t)} + (x_1(t)-x_3(t)-\sin(t))^2 dt +$$ $$+ \int_0^5 \Big( x(t) - \int_0^t z(\tau) d \tau \Big)^2 dt,$$
   where the value $\lambda = 2$ is taken. It is obvious that $z_1^*(t) = 1$, $z_2^*(t)=0$, $z_3^*(t)=1-\cos(t)$, $t \in [0, 5]$, $I(x^*, z^*) = 0$.
   
  Take $(x_{(1)}, z_{(1)}) = (0, 0, 0, 1, 0, 0)' $ as the initial point, then \linebreak $I(x_{(1)}, z_{(1)}) = 44.30267$.  As the iteration number increased, the discretization rank gradually increased during the solution of the auxiliary problem of finding the direction of the subdifferential descent described in the algorithm and in the end the discretization step was equal to $2.5 \times 10^{-2}$. At the 178-st iteration the point 
  $$x_{(178)} = (-0.0000056 t^5+0.000071 t^4-0.000336 t^3+0.000783 t^2+0.999079 t,$$ $$0, 0.00002217 t^9-0.0005342 t^8+0.0050314 t^7-0.023074 t^6+0.053894t^5-$$ $$-0.097232 t^4+0.261281 t^3-0.070276 t^2+0.0352 t)',$$ $$z_{(178)}= (0.000004 t^3-0.000145 t^2+0.000744 t+0.999079,$$ $$0, 0.0001996 t^8-0.0042738 t^7+0.035220 t^6-0.138441 t^5+0.269469 t^4-0.388926 t^3+$$ $$+0.783842 t^2-0.14055 t+0.0352)'$$
   was constructed and the functional value $I(x_{(178)}, z_{(178)}) = 0.0015$. For the convenience of presentation, the Lagrange interpolation polynomial has been given which quite accurately approximates (that is the interpolation error does not affect the value of the functional presented with a given accuracy but (insignificantly) affects the given value of the norm of the smallest subgradient) the resulting trajectory. Herewith, $||\overline{v}_{(178)}||_{L_2^{3}[0,T] \times L_2^{3}[0, T]} \approx 0.0324$.

Since, in fact, $z(t) = \dot x(t)$, $t\in [0,5]$, then one may put  $$x_{(178)} = (-0.0000056 t^5+0.000071 t^4-0.000336 t^3+0.000783 t^2+0.999079 t,$$ $$0, 0.00002217 t^9-0.0005342 t^8+0.0050314 t^7-0.023074 t^6+0.053894t^5-$$ $$-0.097232 t^4+0.261281 t^3-0.070276 t^2+0.0352 t)',$$ $$z_{(178)}= (-0.000028 t^4+0.000284 t^3-0.001009 t^2+0.001565 t+0.999079,$$ $$0, 0.0001996 t^8-0.0042738 t^7+0.035220 t^6-0.138441 t^5+0.269469 t^4-0.388926 t^3+$$ $$+0.783842 t^2-0.14055 t+0.0352)'$$ and then $\overline{J}(x_{(178)}) = I(x_{(178)}, z_{(178)}) \approx 0.00147$, i. e. the error of the functional value does not exceed the magnitude  $10^{-3}$.

 It is interesting to consider the ``behavior'' of the function $z_1(t)$ while increasing the number of iterations. This trajectory is at a ``far'' distance from the true value and insignificantly changes over most of the time interval and noticeably ``improves'' only on a short period of time starting from the finite moment $T=5$. Gradually, the function ``flattens out'' approaching the true value at smaller moments of time and finally becomes close to the solution on the whole time interval. The typical ``behavior'' of the function $z_1(t) - 1$ is illustrated in Figure 2 for some iterations. \newline
 
 \text {\textbf{Table 2}\, Example 4}
$$
\begin{array}{c|c|c|c}
        k&I(x_{(k)}, z_{(k)})&||\overline v(x_{(k)}, z_{(k)})||&||x_{(k)} - x^*||\\
  \hline & &  \\[-3mm]
     1 & 44.3027 & 19.7852 & 9.5249 \\
     50 & 0.8549 & 0.5605 & 2.6881 \\
     100 & 0.0203 & 0.1292 & 0.2327  \\
     150 & 0.0033 & 0.0401 & 0.0201 \\
      178 & 0.0015 & 0.0324 & 0.0189
\end{array}
$$

 \begin{figure*}[h!]
\begin{minipage}[h]{0.3\linewidth}
\center{\includegraphics[width=1\linewidth]{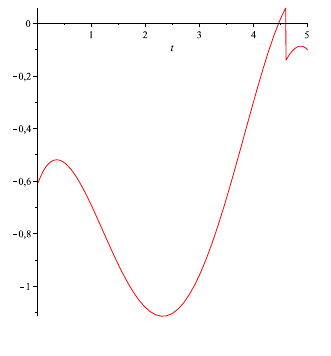} }
\end{minipage}
\hfill
\begin{minipage}[h]{0.3\linewidth}
\center{\includegraphics[width=1\linewidth]{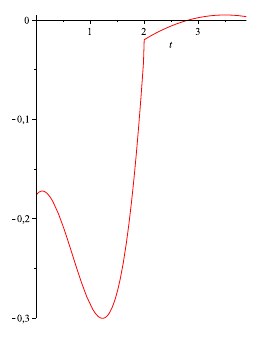} }
\end{minipage}
\hfill
\begin{minipage}[h]{0.3\linewidth}
\center{\includegraphics[width=1\linewidth]{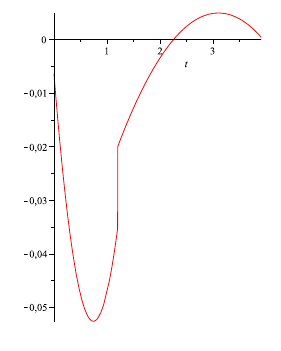} }
\end{minipage}
\caption{Example 4, function: $z_1(t) - 1$, iterations: 12, 68, 87}
\label{ris:image2}
\end{figure*}
 
  {\bf Remark 5.}   
   The examples considered show that in some cases it is required to take a sufficiently large value of the penalty parameter $\lambda$, which can lead to additional computational difficulties. In order to overcome this problem, in future investigations it is planned to consider the functional $$\displaystyle \overline{\varphi}(x,z) =  \int_0^T \Big| x(t) - x_0 -  \int_0^t z(\tau) d \tau \Big| d t$$ instead of the functional $\varphi(x,z)$ which will improve the accuracy of fulfillment of this constraint and will significantly decrease the value $\lambda$. For example, let one consider the functional $\lambda_1 \varphi (x,z)$ (see formula 9)  and take $\varepsilon_1$ as the acceptable error, that is $\lambda_1 \varphi (x,z) = \varepsilon_1$. Then if we put $\lambda_2 \overline\varphi (x,z) = \varepsilon_2$, then with the help of H\"{o}lder's inequality it is easy to check that the penalty parameters are related by the following relation $\displaystyle{\lambda_1 = \frac{\varepsilon_1 T^2 \lambda_2^2}{2 \varepsilon_2^2}}$, so in most practical cases the value of $\lambda$ can be reduced by several orders of magnitude. (Take, e.g., $\lambda_1 = 100$, $\varepsilon_1 = \varepsilon_2 = 10^{-2}$, $T=1$, then $\lambda_2 = \sqrt{2}$.) The hypothesis is that the functional $I(x, z)$ (with the functional $\overline{\varphi}(x,z)$) may turn out to be an exact penalty one \cite{DolgExact}, which means (speaking not strictly and without formal definitions) the possibility to take even the finite value of the parameter $\lambda^*$ such that for all $\lambda > \lambda^*$ the minimizer of this functional will strictly satisfy the corresponding constraint.


\section{Acknowledgements}
The author is sincerely grateful to his colleagues Maksim Dolgopolik and \linebreak Grigoriy Tamasyan for numerous fruitful discussions and to the anonymous referee, whose comments helped to significantly improve the paper.

The work was supported by the Russian Science Foundation (project no. 21-71-00021).





\end{document}